\documentclass[11pt]{amsart}

\pdfoutput=1

\usepackage[text={420pt,650pt},centering]{geometry}

\usepackage{cancel}
\usepackage{esint,amssymb} 
\usepackage{graphicx}
\usepackage{MnSymbol}
\usepackage{mathtools} 
\usepackage[colorlinks=true, pdfstartview=FitV, linkcolor=blue, citecolor=blue, urlcolor=blue,pagebackref=false]{hyperref}
\usepackage{microtype}

\usepackage{bm}
\usepackage{dsfont}
\usepackage{mathrsfs}
\usepackage{xcolor}

\usepackage{tcolorbox}
\tcbuselibrary{most}
\usepackage{wrapfig}
\newcounter{myboxcounter}
\definecolor{lightblue}{RGB}{173, 216, 230}

\newtcolorbox[auto counter, number within=section]{mybox}[1][]{
    colback=lightblue!20,   
    colframe=lightblue!20,  
    width=\textwidth,       
    boxrule=0mm,            
    arc=0mm,                
    boxsep=3mm,             
    fonttitle=\small\bfseries, 
    fontupper=\small,       
    before upper={\refstepcounter{myboxcounter}\textbf{Box~\themyboxcounter.} \textbf{#1}\\}, 
    breakable,              
}

\usepackage{caption}

\newtheorem{proposition}{Proposition}
\newtheorem{theorem}[proposition]{Theorem}

\newtheorem{conjecture}[proposition]{Conjecture}
\theoremstyle{remark}

\theoremstyle{definition}

\numberwithin{equation}{section}
\numberwithin{proposition}{section}
\numberwithin{figure}{section}
\numberwithin{table}{section}
\numberwithin{myboxcounter}{section}

\newcommand{\R}{\mathbb{R}}

\newcommand{\E}{\mathbb{E}}

\newcommand{\ep}{\varepsilon}
\newcommand{\eps}{\varepsilon}

\renewcommand{\le}{\leqslant}
\renewcommand{\ge}{\geqslant}

\renewcommand{\bar}{\overline}

\newcommand{\Ll}{\left}
\newcommand{\Rr}{\right}
\renewcommand{\d}{\mathrm{d}}
\newcommand{\dr}{\partial}

\newcommand{\1}{\mathbf{1}}

\newcommand{\mcl}{\mathcal}

\newcommand{\de}{\delta}

\newcommand{\la}{\left\langle}
\newcommand{\ra}{\right\rangle}

\newcommand{\bip}{{\mathrm{bip}}}

\begin{document}

\author{Jean-Christophe Mourrat}
\address[Jean-Christophe Mourrat]{Department of Mathematics, ENS Lyon and CNRS, Lyon, France}

\subjclass[2010]{}
\date{\today}

\title[An informal introduction to the Parisi formula]{An informal introduction to the Parisi formula}

\begin{abstract}
This note is an informal presentation of spin glasses and of the Parisi formula. We also discuss some models for which the Parisi formula is not well-understood, and some partial progress that relies upon a connection with partial differential equations. 
\end{abstract}

\maketitle

Statistical mechanics aims to model the emergent properties of systems that are made of a large number of elements. One can think of a gas made of many small particles, which gives rise to macroscopic concepts such as the density, the pressure, or the temperature of the gas. 
Here we will focus on a class of models of statistical mechanics called \emph{spin glasses}. The critical feature of these models is that there is a lot of ``disagreement'' between the elementary units of the system, as these models aim to capture features of ``complex'' systems. Spin glasses have inspired many developments in a variety of topics including statistics, computer science, high-dimensional geometry, or combinatorics, and many facets of these models have been studied. Our main goal here is to discuss some of the ideas surrounding a fundamental result called the Parisi formula, which identifies the limit free energy of some of these models. The free energy is a natural quantity from a physicist's perspective, and can be thought of as a Laplace transform of the variables of interest, so the identification of its limit yields rich insight into the behavior of the model. We will also discuss an intriguing connection between the Parisi formula and certain partial differential equations, and how this may help to address some open problems. 

The outline of this note is as follows. In the first section, we introduce one of the most basic models of a spin glass called the Sherrington-Kirkpatrick model, and explain what features of this model make it indeed a spin glass. Section~\ref{s.parisi} presents the Parisi formula per se; we also discuss the crucial insight that the support of the associated Gibbs measure is approximately ultrametric. In Section~\ref{s.general}, we discuss extensions to more general models, and showcase one example called the bipartite model that is currently less well-understood. In Section~\ref{s.pdes}, we sketch how to rephrase the Parisi formula using certain partial differential equations, as well as some partial results concerning the bipartite model.

This note has a number of footnotes and of pieces of text that are within colored boxes. These are asides or more technical discussions that you can freely skip if you so wish. Reading up until the end of Section~\ref{s.parisi} should give you a good idea of what the Parisi formula is about, and you may decide to stop there. The remaining sections relate more directly to my recent research activity.

%
%
%
%
%
%
\section{The Sherrington-Kirkpatrick model}

To get a clearer sense of what spin glasses are, let us introduce a paradigmatic example called the Sherrington-Kirkpatrick (SK) model  \cite{sherrington1975solvable}. 
It may help to imagine ourselves being confronted with the following situation. Suppose that we have $N$ individuals $\{1,\ldots, N\}$ that we need to split into two groups. We can represent an assignment into two groups as a vector $\sigma \in \{\pm 1\}^N := \{-1,1\}^N$, where $\sigma_i$ indicates to which group the individual indexed by $i$ is assigned. For each pair $(i,j)$, we are given a number $W_{ij}$ which represents how much the individual $i$ likes the individual $j$; it is positive and very large if the individual $i$ really likes the individual $j$, while it is very negative if the individual $i$ really dislikes the individual $j$. We want to look for an assignment $\sigma \in \{\pm 1\}^N$ that makes the following quantity as large as possible:
\begin{equation}  
\label{e.def.HN}
H_N(\sigma) := \frac 1 {\sqrt{N}} \sum_{i,j = 1}^N W_{ij} \sigma_i \sigma_j.
\end{equation} 
It might have been more natural to write $\1_{\{\sigma_i = \sigma_j\}}$ in place of $\sigma_i \sigma_j$, but since $\1_{\{\sigma_i = \sigma_j\}} = \frac{1}{2}(\sigma_i \sigma_j+1)$, this is inconsequential\footnote{One reason for defining the model using $\sigma_i \sigma_j$ in place of $\1_{\{\sigma_i = \sigma_j\}}$ has to do with other motivations coming from the modelling of magnetic alloys. From a purely mathematical perspective, it is also more pleasant to think of $H_N$ as a polynomial function of $\sigma$, which we can think of as being defined everywhere in $\R^N$ instead of just on $\{\pm 1\}^N$. Another small point is that I find it slightly more convenient not to assume that $W_{ij}$ equals $W_{ji}$, but this is also a detail.}. The normalization $\frac 1 {\sqrt{N}}$ will be convenient later on.

\begin{wrapfigure}{r}{0.3\linewidth}
    \centering
    \includegraphics[scale=1]{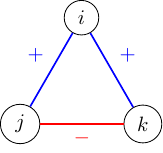}
    \captionsetup{width=0.95\linewidth}
\caption{{\small A simple si\-tuation with frustration. Here the coefficients $(W_{ij})$ suggest to set $\sigma_i = \sigma_j$, $\sigma_i = \sigma_k$, and $\sigma_j = - \sigma_k$, but we cannot realize these three conditions simultaneously.}}
    \label{f.frustration}
\end{wrapfigure}

One quickly realizes that finding the configuration $\sigma \in \{\pm 1\}^N$ that realizes the maximum of $H_N$ is not going to be straightforward. Even when only three individuals $i$, $j$ and $k$ are involved, we may be in a situation as depicted in Figure~\ref{f.frustration}: maybe $W_{ij}+W_{ji} > 0$ and $W_{ik}+W_{ki} > 0$, which would suggest to assign $i$, $j$ and $k$ to the same group; but if $W_{jk} + W_{kj} < 0$, then the individuals $j$ and~$k$ would rather be in different groups, and there is no way to reconcile each of the pairwise preferences. 

In other words, certain pairs will typically end up being  \emph{frustrated}. In order to find the optimal configuration, some compromises need to be made, and a close inspection of the coefficients $(W_{ij})$ is required. 
More generally, for $N$ large, naive quick strategies\footnote{For example, as long as there exists an index \( i \) such that changing \( \sigma_i \) to \( -\sigma_i \) increases \( H_N \), we look for the index that produces the greatest difference, and we iterate.} for finding a configuration $\sigma$ such that $H_N(\sigma)$ is large will typically get stuck at a local maximum of the mapping $H_N$, and will fail to reach the true maximizer. The presence of these frustrations, and the related fact that simple methods typically do not succeed in finding the maximizing configuration, are the defining features of \emph{glassy} systems. How this relates to the glass of our daily lives is discussed in Box~\ref{b.glass}. The variables $\sigma_1,\ldots, \sigma_N$ are often called \emph{spins}, because much of statistical mechanics has focused on the modeling of magnetic materials, and the model is therefore called a \emph{spin glass}. 

\begin{mybox}[Why the word ``glass'' for these models?]
\label{b.glass}
In order to create glass, one starts by heating up silica (i.e. sand) and a bit of calcium and sodium carbonates (i.e.\ lime and soda) until they melt. One key aspect of the fabrication of glass is that the liquid mixture then needs to be cooled rapidly (``quenched'' is the technical term), so as to ``trap'' the microscopic configuration into a disordered state inherited from the liquid phase. After the rapid quench, the material is forever trying to slowly evolve towards its energetically-preferred organized state, but encounters locally frustrated configurations that are increasingly difficult to overcome. So in this case the frustrations and slow dynamics emerge from the intricate geometry of the configuration of particles. The conflicting coefficients $(W_{ij})$'s in the SK model make it much more analytically tractable than any model that would really try to be faithful to the ``geometric disorder'' of a true glass, but one hopes that the two settings share certain qualitative properties. 
\end{mybox}

In line with the necessary presence of these frustrations, one can show that  the problem, given the coefficients $(W_{ij})$, of finding a configuration $\sigma \in \{\pm 1\}^N$ that maximizes $H_N$, is NP-hard in general. In fact, the problem is NP-hard even if we only aim to find a configuration $\sigma \in \{\pm 1\}^N$ such that $H_N(\sigma)$ is at least a fixed positive fraction of the maximal value, no matter how small we allow the ratio to be~\cite{arora2005non}. We will however depart from this worst-case sort of analysis, by focusing instead on what a ``typical'' instance of the problem looks like. There are surely different possible ways to clarify what ``typical'' means here, but we choose to postulate that the coefficients $(W_{ij})_{i,j \le N}$ are independent Gaussian random variables with mean zero and unit variance. The key hypothesis here is that these are independent and all have the same mean and variance; the hypothesis that they are Gaussian is a convenience but would not change the fundamental properties that will be discussed below. 


To sum up, we let $(W_{ij})_{i,j \ge 1}$ be independent centered Gaussian random variables with unit variance, we define $H_N$ according to \eqref{e.def.HN}, and we aim to study quantities such as
\begin{equation}
\label{e.def.max}
\frac 1 N \max_{\sigma \in \{\pm 1 \}^N} H_N(\sigma),
\end{equation}
in the limit of large $N$. More generally, we would be interested in understanding the geometry of the function~$H_N$, for instance with a view towards finding configurations~$\sigma$ that essentially realize the maximum in \eqref{e.def.max}. A particularly fruitful way to probe this is to consider a family of probability measures associated with $H_N$ called \emph{Gibbs measures}; see also Box~\ref{b.gibbs}. In our context, given a parameter $\beta \ge 0$, the Gibbs measure at ``inverse temperature''~$\beta$ is the probability measure that attributes to each $\sigma \in \{\pm 1\}^N$ a probability proportional to $\exp(\beta H_N(\sigma))$.\footnote{Physicists prefer to add a minus sign here, writing $\exp(-\beta H_N(\sigma))$ in place of $\exp(\beta H_N(\sigma))$, but since the laws of $H_N$ and $-H_N$ are identical, this does not really matter and I prefer to avoid the proliferation of minus signs. The slight drawback of this convention is that the system now has a preference for \emph{large} values of $H_N$, while physicists prefer to think that the energy function $H_N$ ought to be minimized.} For each value of $\beta \ge 0$, this probability measure essentially concentrates on a level set of $H_N$, in the sense that $H_N/N$ is essentially constant under the Gibbs measure; and as $\beta$ is taken larger and larger, the measure becomes concentrated on near-maximizers of $H_N$. In order to understand this family of Gibbs measures, it is very fruitful to elucidate the behavior of the quantity
\begin{equation}
\label{e.def.free.energy}
F_N(\beta) := \frac 1 N \E \log \sum_{\sigma \in \{\pm 1\}^N} \exp(\beta H_N(\sigma)).
\end{equation}
The expectation $\E$ is with respect to the randomness coming from the coefficients $(W_{ij})$.
The quantity $F_N(\beta)$ is usually called the \emph{free energy} of the system\footnote{even though physicists would rather divide $F_N$ by $\beta$ and add a minus sign (to match the one they would have in the exponential) before calling it that.}. Understanding this quantity or generalizations of it is very rich in insight concerning the associated Gibbs measure\footnote{The Gibbs measure depends itself on the sampling of the random coefficients $(W_{ij})$, and we hope to describe typical properties with respect to this sorting, or average properties.}. For mathematicians, this makes intuitive sense if we think of it as a sort of log-Laplace transform of the function $H_N$; we can also see for instance that the derivative in $\beta$ of $F_N$ gives us access to the average of $H_N$ under the Gibbs measure, averaged also over the coefficients $(W_{ij})$. We can alternatively think of $F_N(\beta)$ as a soft version of the maximum in \eqref{e.def.max}, so that this maximum is approximately $F_N(\beta)/\beta$ if we take $\beta$ sufficiently large\footnote{The key step to justify this point is to show that this maximum is close to its expectation; see for example \cite[Exercises 6.2 and 6.3]{HJbook}.}.

\begin{mybox}[Gibbs measures]
\label{b.gibbs}
We may as well just think of Gibbs measures as useful mathematical tools to probe the level sets of $H_N$, but when $H_N$ represents the energy function of a real physical system, one can argue from first principles that at equilibrium in an environment at inverse temperature $\beta$, the system will indeed be distributed according to this measure. In fact, this is a general phenomenon that can also show up in purely mathematical contexts. For intance, suppose that there are variables $x_1,\ldots, x_N$ taking values in a finite set $\{e_1, \ldots, e_K\}$, say with $e_1 < \ldots < e_K$, and let $\bar e \in (e_1, e_K)$. How does one of the variables, say $x_1$, look like if we pick the whole vector $(x_1,\ldots, x_N)$ uniformly at random among all those that satisfy $N^{-1} \sum_{i = 1}^N x_i \simeq \bar e$ ? (Here we write $\simeq$ to allow for some wiggle room so as to make sure that there exist such vectors, for instance we can ask that the difference between the terms on the two sides of $\simeq$ is at most~$\eps_N$, with $\eps_N \to 0$ and $N \eps_N \to +\infty$ as $N$ tends to infinity.) One can show that as $N$ tends to infinity, the probability that $x_1 = e_k$ is proportional to $\exp(-\beta e_k)$, where $\beta$ is such that 
\begin{equation}  
\label{e.def.beta}
\frac{\sum_{k = 1}^K e_k \exp(-\beta e_k)}{\sum_{k = 1}^K  \exp(-\beta e_k)}  = \bar e. 
\end{equation}
In other words, asymptotically as $N$ tends to infinity, the law of $x_1$ converges weakly to the Gibbs measure at inverse temperature $\beta$, where $\beta$ is defined by \eqref{e.def.beta}. 
One may for instance consult \cite[Section~1.1]{HJbook} for more on this and on the relationship between the Gibbs measure, the free energy, and the entropy of the Gibbs measure.
\end{mybox}

It is not very difficult to show that the quantity in \eqref{e.def.max}, and then also that in \eqref{e.def.free.energy}, stay away from zero and infinity as $N$ tends to infinity\footnote{You can try! Here are some hints. For the upper bound, observe that for each $\sigma \in \{\pm 1\}^N$, the random variable $H_N(\sigma)$ is a centered Gaussian with variance $N$, so  we can easily estimate the probability that $H_N(\sigma)$ is above $CN$ for each $\sigma$ individually. For the lower bound, one can rely on a naive greedy procedure in which, assuming we have already commited to a choice of $\sigma_1,\ldots, \sigma_i$, we pick the choice of $\sigma_{i+1}$ that maximizes that part of $H_N$ that only involves $\sigma_1,\ldots, \sigma_{i+1}$, and we iterate over $i$.}. Can we determine the limits of these quantities? Despite a somewhat lengthy introduction, I hope that you can appreciate the simplicity of the question. And yet, something I find fascinating is that the answer to this question is incredibly rich and complex.

%
%
%
%
%
%
\section{The Parisi formula}
\label{s.parisi}

An initial guess for the limit of the free energy $F_N(\beta)$ in \eqref{e.def.free.energy} was proposed by physicists in the paper that introduced the model \cite{sherrington1975solvable}, but it was already understood there that the proposed answer could not be valid for large values of $\beta$, i.e.\ at low temperature. In one of his most celebrated contributions, Giorgio Parisi then came up with a sophisticated non-rigorous procedure, called the replica method, that led to what is now called the \emph{Parisi formula} \cite{parisi1979infinite, parisi1980order, parisi1980sequence, parisi1983order}. This formula is described in full in Box~\ref{b.parisi}; here we content ourselves with the fact that it takes the form
\begin{equation}
\label{e.parisi.simple}
\lim_{N \to \infty} F_N(\beta) = \inf_{\mu \in \Pr([0,1])} \mcl P(\mu),
\end{equation}
for some functional $\mcl P$, where $\Pr([0,1])$ denotes the space of probability measures on~$[0,1]$. 

\begin{mybox}[The Parisi formula in full]
\label{b.parisi}
The Parisi formula states that
\begin{equation}\label{e.SK.Parisi}
\lim_{N\to+\infty}F_N(\beta)=\inf_{\mu \in \Pr([0,1])}\bigg( \Phi_\mu(0,0) -\beta^2\int_0^1 t\mu([0,t])\d t+\log(2)\bigg),
\end{equation}
where $\Phi_\mu :[0,1]\times \R \to \R$ is the solution to the backwards parabolic equation
\begin{equation}\label{e.SK.Parisi.PDE.dis}
\begin{cases}
-\partial_t\Phi_\mu(t,x)=\beta^2\Big(\partial_x^2\Phi_\mu(t,x)+\mu([0,t])\big(\partial_x\Phi_\mu(t,x)\big)^2\Big) & \text{for } (t,x) \in [0,1]\times \R,\\
\Phi_\mu(1,x)=\log\cosh(x)& \text{for } x\in \R.
\end{cases}
\end{equation}
\end{mybox}

There are many shades of non-rigorous arguments. Some do not have all the $\eps$'s and~$\delta$'s or discard some annoying but plausibly negligible terms, yet most mathematicians would feel rather convinced by them. Giorgio Parisi's techniques were of a different sort though, see Box \ref{b.replica}. In fact, the opinions on the validity of his prediction initially varied in the physics community. Besides the rather creative nature of the arguments involved, part of the initial skepticism may have come from the fact that the Parisi formula is a minimization problem, as displayed in \eqref{e.parisi.simple}. For each fixed $N$, the free energy can be rewritten as a supremum over probability measures of an ``energy'' term and an ``entropy'' term\footnote{While we did not proceed to fully explain this in Box~\ref{b.replica}, the free energy and the entropy are in a convex duality relationship; see also \cite[(1.14)]{HJbook} and \cite[Corollaries~4.14 and 4.15]{blm} for more on this.}; once one gets familiar with it, this formulation feels very intuitive and appealing. One could then expect that the task of identifying the limit of the free energy boils down to understanding how to simplify this representation in the limit of large $N$, while preserving the main structure as a supremum. And yet the formula in \eqref{e.parisi.simple} displays an infimum instead\footnote{It does not help that in the replica calculation as explained in Box~\ref{b.replica}, a maximization problem appears that is then flipped into a minimization problem as the positive integer $n$ becomes smaller than~$1$, with not much of an explanation besides the fact that one would otherwise end up with something nonsensical. I cannot help pointing out that a rewriting of the Parisi formula that takes the form of a supremum was recently found in \cite{mourrat2024uninverting}.}.

\begin{mybox}[The replica method]
\label{b.replica}
The replica trick aims to exploit the fact that $\log x = \lim_{n \to 0} (x^n-1)/n$ together with the calculation of what in our context would read as
\begin{equation}  
\label{e.def.moments}
\E\bigg[\bigg(\sum_{\sigma \in \{\pm 1\}^N} \exp(\beta H_N(\sigma))\bigg)^n\bigg],
\end{equation}
where $n$ is a positive integer. The advantage of working with an integer $n$ is that we can then expand the power and rewrite the expression in \eqref{e.def.moments} as a sum over $n$ copies of the variable $\sigma$ of some expected value we can compute; the $n$ variables ranging in $\{\pm 1\}^N$ are usually called ``replicas'' (although we will use this word in a slightly different sense further below). This replica trick was already used in the paper that introduced the model \cite{sherrington1975solvable} and earlier for other models; a short survey of it is in \cite[Appendix]{parisi2023nobel}. Giorgio Parisi created the art of sending $n$ to zero in just the right way as one goes along the calculation so as to arrive at the correct answer for the limit of the free energy. To give a flavor of the manipulations involved, here is a quote from his Nobel lecture \cite{parisi2023nobel}. ``After many trials I had an intuition: in other papers the $n$ indices were
divided into $n/m$ groups of $m$ elements each [...]. Everybody
was assuming that $m$ was an integer [...].
I made the bold assumption that $m$ could be a noninteger
number, more precisely a number in the interval $[0,1]$. For
example, I was dividing the $n$ replicas into $2n$ groups of $1/2$
replicas each. Of course, that is crazy, but my viewpoint was I
should first check if this crazy idea was leading to correct
results and postpone other questions to a later stage. [...] At the end of the paper, I added the observation that one
could improve the theory by dividing the $n/m$ groups into
$(n/m)/m_1$ groups of replicas, where $m_1$ was a new variational
parameter. I was also conjecturing the correct solution was
obtained when the procedure was repeated an infinite number
of times. Some fancy group theory arguments also were added, arguing the permutation group of zero objects is an infinite group because it contains itself as a proper subgroup. [...] The response of the referee was remarkable. In a nutshell:
\emph{The approach does not make sense, but the numbers coming
from the formulae are reasonable, so it can be published. The last observation is not worth the paper on which it is
written and it should be removed.} I laughed because in the meanwhile I extended the computation to an infinite number of subdivisions''. That latter calculation led to what we now call the Parisi formula.
\end{mybox}

\begin{figure}[h]
\centering
\includegraphics[width=0.55\linewidth]{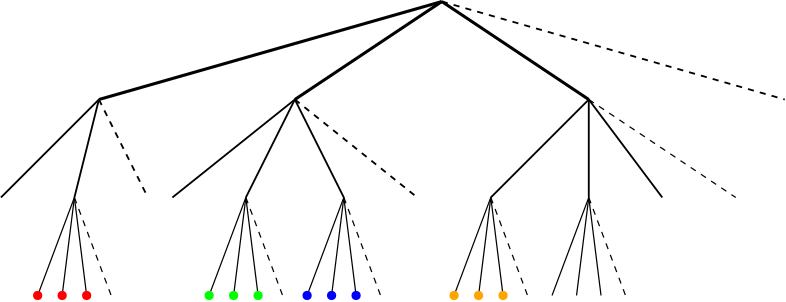}
\captionsetup{width=0.95\linewidth}
\caption{{\small Except for some choices of $(W_{ij})$ of small probability and for large $N$, samples from the Gibbs measure essentially behave as if they were sampled from an ultrametric space. An ultrametric space can be encoded on the leaves of a tree as in the picture above, where the blue points are at the same distance from one another; they are themselves all at the same distance from any green point; the blue, green and orange points are all at the same distance from any red point, etc.}}
\label{f.ultrametric}
\end{figure}

Further work by physicists later elucidated a wealth of new information that was consistent with the formula in \eqref{e.parisi.simple} \cite{mezard1984replica,mezard1984nature, mpvbook}. They discovered that, as $N$ tends to infinity (and except for some choices of $(W_{ij})$ of small probability), the associated Gibbs measure has a very complex and yet also very precise hierarchical structure, and that the minimizing measure in \eqref{e.parisi.simple} fully describes it. One key aspect of this description is that as $N$ becomes very large, the Gibbs measure is essentially supported on an ultrametric set; see Figure~\ref{f.ultrametric} and Box~\ref{b.ultrametricity} for more on this. 

\begin{mybox}[Ultrametricity of the Gibbs measure]
\label{b.ultrametricity}
We wish to discuss the behavior of a probability measure (the Gibbs measure) defined on $\{\pm 1\}^N$, or more generally on $\R^N$, in the limit of large $N$. Since the space over which this probability measure is defined changes with $N$, one cannot simply ask whether the sequence converges weakly for instance. There is also a lot of invariance in the model, such as by permutations of the coordinates, so we may want to discard these symmetries while keeping interesting information. In a sense, we aim to retain information about the Gibbs measure as an abstract random metric space, ignoring how it is embedded into $\R^N$ and discarding small-probability parts. One very convenient way to do this is to study the convergence of the joint law of the relative distances between finitely many independent samples from the measure\footnote{This idea defines a topology on the space of (isometry classes of) random measure spaces called the Gromov-weak topology, and this topology is equivalent to that induced by the Gromov-Prokhorov metric \cite{greven2009convergence}.}; these independent samples are called \emph{replicas}. 

Recall that, for each choice of the parameters $(W_{ij})$ and each $\beta \ge 0$, the Gibbs measure at inverse temperature $\beta$ is the probability measure that attributes to the configuration $\sigma \in \{\pm 1\}^N$ a weight proportional to $\exp(\beta H_N(\sigma))$. It is thus a \emph{random} probability measure, as it depends on our choice of the coefficients $(W_{ij})$, and we are going to discuss properties of this probability measure that occur with high probability over the choice of these coefficients $(W_{ij})$ and for large $N$. What physicists discovered, and was later proved rigorously (up to allowing ourselves a small perturbation of the energy function that does not change the value of the limit free energy), is that the Gibbs measure is asymptotically \emph{ultrametric}. In Euclidean space, a set $S$ is ultrametric if it satisfies the following stronger form of the triangle inequality: for every $x,y,z \in S$, we have 
\begin{equation*}  
|x-y| \le \max(|x-z|, |y-z|).
\end{equation*}
This can be rephrased as the fact that whenever we draw a triangle with endpoints in $S$, the two largest sides of the triangle have the same length. 
In our context, the approximate ultrametricity of the Gibbs measure means that the following holds for every $\eps > 0$ and $\beta \ge 0$ with probability tending to $1$ in the $W$'s as $N$ tends to infinity: if we pick three independent samples (``replicas'') $\sigma$, $\sigma'$, and $\sigma''$ from the Gibbs measure at inverse temperature $\beta$, then with Gibbs probability tending to $1$ as $N$ tends to infinity, we have
\begin{equation*}  
|\sigma - \sigma'| \le \max(|\sigma - \sigma''|, |\sigma' - \sigma''|) + \eps\sqrt{N}
\end{equation*}
(where $|\cdot|$ is the Euclidean distance in $\R^N$, and there is a natural $\sqrt{N}$ scaling factor since $|\sigma| = \sqrt{N}$).
Any ultrametric set can be represented as the leaves of a tree, with the distance between two leaves being only allowed to depend on the depth of the most recent common ancestor of the two leaves; see Figure~\ref{f.ultrametric} (or e.g.\ \cite[Lemma~4.2]{chen2024free}). The Gibbs measure thus splits into a number of ``pure states''\footnote{This is vaguely analogous to the decomposition of a probability measure into its ergodic components in the theory of dynamical systems.} represented by the leaves of the tree; any two configurations sampled from a pure state are at a fixed distance from one another (with high probability). These pure states organize themselves into clusters, so that the distance between two configurations sampled from two different pure states in the same cluster is essentially constant; these clusters organize themselves into super-clusters; and so on, over a potentially infinite hierarchy of clusters and super-clusters. The exact number of hierarchies depends on the specifics of the problem. There will be no such hierarchy at high temperature (low $\beta$), in some sense the tree collapses into a single point in this case. At low temperature (high $\beta$), physicists seem to expect that most ``natural'' problems display either one or an infinite number of levels of hierarchy, although we know how to engineer spin-glass models with an arbitrary number of them.

Beyond the fact that the asymptotic structure is ultrametric, in fact the full metric structure as well as the probability weights that the Gibbs measure attributes to each leaf are completely characterized in terms of the minimizer in \eqref{e.parisi.simple}. It would be too long to describe this structure precisely here (see for instance \cite{pan} or \cite{HJbook}, including \cite[Section~6.3]{HJbook} which explains how a 1-level hierarchy emerges from a much simpler toy model). I do want to stress though that I find this rigidity really remarkable. In particular, up to some technical caveats, we have that the hierarchical structure has depth $K$ if and only if the support of the measure $\mu$ that minimizes the Parisi formula in~\eqref{e.parisi.simple} contains exactly $K+1$ points. Moreover, if we pick two independent samples $\sigma$ and $\sigma'$ from the Gibbs measure, then the law of their ``overlap'' $\sigma \cdot \sigma'/N$ converges weakly to $\mu$. 

Whenever a non-trivial hierarchical structure appears, physicists say that there is \emph{replica-symmetry breaking}; and if there are $K$ levels to the hierarchy, they would say that the system has $K$ levels of replica-symmetry breaking. 

To explain the choice of language, it is useful to ask ourselves first what physicists mean when they say that a symmetry is broken. Suppose that I try to balance a pencil that I place vertically on the tip of my finger. By symmetry, the only thing that can happen is that the pencil stays perfectly vertical, right? Well, my experiments do not match this prediction, as the pencil actually quickly falls in some direction. So one could say that the rotation symmetry is broken. Of course, if I try to be a bit diligent in my attempts at balancing the pencil, then the direction in which the pencil falls will be essentially uniformly random among all possible directions. The symmetry is therefore not really broken in a mathematical sense. The correct mathematical model of the outcome is probabilistic, but we have the impression that the symmetry is broken if we ignore this and only do one experiment. 

Coming back to spin glasses and their ultrametric structure: suppose that we pick three independent samples $\sigma$, $\sigma'$ and $\sigma''$ according to the Gibbs measure. Of course there is perfect permutation symmetry between these samples, so surely the distance between $\sigma$ and $\sigma'$ should be the same as that between $\sigma$ and $\sigma''$? Well, this is in fact not the case (see Figure~\ref{f.ultrametric}), so the permutation symmetry between the replicas is broken. The key surprising thing here is that these relative distances remain random even in the limit of large $N$, even though one may have thought otherwise at first, because many quantities become essentially deterministic (``self-averaging'' as physicists would say) in the high-dimensional regime.
\end{mybox}

This series of discoveries generated a lot of excitement; people had identified a toy model of a ``complex'' system, with a very rich and rugged energy landscape, that they could study with great precision using analytical methods. Yet, says Giorgio Parisi \cite{parisi2023nobel}, ``it was possible that the correct [solution] was different and more complex [...]. It was difficult to conclude in a definite way.'' The validity of the Parisi formula then became a certainty with a series of rigorous mathematical contributions, starting with Francesco Guerra \cite{gue03} who proved that $\limsup_{N \to \infty} F_N(\beta) \le \inf_{\mu \in \Pr([0,1])} \mcl P(\mu)$. 
In a mathematical tour de force, Michel Talagrand \cite{Tpaper} then managed to prove the converse bound and therefore give a complete justification to Parisi's formula\footnote{This is one of the key results celebrated in the citation for Michel Talagrand's 2024 Abel prize.}. An alternative proof that covers a broader class of models was later developed by Dmitry Panchenko \cite{pan.aom, pan}. Besides its greater generality, this alternative proof is also more conceptual, and its key step is very interesting on its own, as it consists in justifying the ultrametricity of the underlying Gibbs measure (up to a small perturbation of the energy function).

%
%
%
%
%
%
\section{Towards more general models}
\label{s.general}

Thanks to the insights coming from the Parisi formula and its proof, as well as further developments, we now understand many aspects of the Sherrington-Kirkpatrick model and the structure of its Gibbs measures. To give just one example, there has been much recent progress on whether one can find a polynomial-time algorithm that identifies a configuration $\sigma \in \{\pm 1\}^N$ that essentially maximizes the function $H_N$ (see \cite{ams, gamarnik2021overlap-survey, gamarnik2021overlap-paper, gamarnik2022disordered, huang2021tight, jekel2024potential, montanari2021optimization,  sellke2021optimizing, sub18})\footnote{Roughly speaking, the general criterion is that this is possible if and only if the support of the minimizer $\mu \in \Pr([0,1])$ in the Parisi formula \eqref{e.parisi.simple} is connected for all $\beta$ sufficiently large.}.

The insights that were gained for the SK model have allowed us to make progress on a large variety of other ``complex'' problems that share similar features with ``frustrations'', from statistics and high-dimensional geometry to computer science and combinatorics. Examples include random constraint satisfaction problems \cite{ding2015proof, ding2016satisfiability,  krzakala2007gibbs, mezmon, mezard2002analytic, monasson1999determining}, the random assignment and traveling salesman problems \cite{aldous1, aldous2, mezard1986replica},
community detection and more general problems of large-scale statistical learning \cite{abbe2017community, zdekrz} (see also \cite[Section~4]{HJbook}), error correcting codes in information theory \cite{richardson2008modern}, and combinatorial problems such as graph coloring \cite{coja2018, ding2016maximum, mulet2002coloring}.

\begin{wrapfigure}{r}{0.3\textwidth}
\centering


\includegraphics[width=0.25\textwidth]{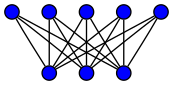}
\captionsetup{width=0.95\linewidth}
\caption{{\small Elementary units are organized into two layers and only interact across layers.}}
    \label{f.bipartite}
\end{wrapfigure}

In this presentation, I will single out a seemingly modest generalization of the SK model in which the variables are of two different types. We can visualize this by thinking of the variables as being organized into two layers as depicted in Figure~\ref{f.bipartite}, with each layer being made of variables of a single type; we now postulate that there are direct interactions between variables of different types only. I find this model interesting because it is closely related to several classical models of artificial neural networks. One is the Hopfield model, which is a model of memory storage and retrieval with a long history \cite{amari1972learning, amit1985spin, amit1987statistical, barra2012equivalence, hopfield1982neural, little1974existence, Tbook2}. Another one is the so-called restricted Boltzmann machine, which is an artificial neural-network architecture that was popular until around ten years ago for learning data distributions and then generating new samples \cite{hinton2012practical, smolensky1986information}; in this case our model is only capturing the initial stage of the network before any learning occurs, but physicists have already managed to gain very interesting insight from there \cite{barra2018phase, tubiana2018restricted, tubiana2017emergence}\footnote{In the time since this paragraph was written, the Nobel committee awarded the 2024 Physics prize to John Hopfield and Geoffrey Hinton, citing in particular their works on these models.}.

To formalize the model precisely, we can represent the variables as a pair $\sigma = (\sigma_1, \sigma_2) = (\sigma_{1,1},\ldots, \sigma_{1,N}, \sigma_{2,1}, \ldots, \sigma_{2,N}) \in \{\pm 1\}^N \times \{\pm 1\}^N$, and we set
\begin{equation}  
\label{e.def.HN.bip}
H_N^\mathrm{bip}(\sigma) := \frac 1 {\sqrt{N}} \sum_{i,j = 1}^N W_{i,j} \sigma_{1,i} \sigma_{2,j}.
\end{equation}
For simplicity we have imposed $\sigma_1$ and $\sigma_2$ to be vectors of the same length $N$, but this is not fundamental; the important point is to make sure that the respective sizes of the two layers remain proportional to one another as we consider larger and larger systems. We have also retained the idea that each of the variables takes values in $\{\pm 1\}$, but this can also be relaxed\footnote{By sticking with the simplest possible version of the model, we have artificially introduced a symmetry between the two layers, but we want to insist on devising analysis techniques that do not rely on this symmetry; a simple way to break the symmetry while staying within our context and notation is to add a term $h \sum_{i=1}^N\sigma_{1,i}$ into the definition of $H_N^{\mathrm{bip}}(\sigma)$, for some parameter $h \neq 0$.}. We will call this model the \emph{bipartite model}. 

This model may seem barely distinguishable from the SK model at first sight. Yet, to this day, we do not know what the limit of the free energy is in this case; in fact, even the fact that the free energy converges as $N$ tends to infinity is not known. To be clear, the free energy here is 
\begin{equation*}  
\frac 1 N \E \log \sum_{\sigma \in \{\pm 1\}^N \times \{\pm 1\}^N} \exp(\beta H_N^\mathrm{bip}(\sigma)).
\end{equation*}
The asymptotic behavior of the maximum
\begin{equation*}  
\frac 1 N \max_{\sigma \in \{\pm 1\}^N \times \{\pm 1\}^N} H_N^\mathrm{bip}(\sigma)
\end{equation*}
is also poorly understood. The problem here goes beyond that of fixing some technical part in the proof of the Parisi formula. Indeed, one may imagine several possible ways for the Parisi formula to generalize to this bipartite model, however it can be shown that none of those candidates for the limit are valid \cite[Section~6]{mourrat2020nonconvex}. 

In order to clarify what distinguishes the bipartite model from the SK model at the technical level, it is best to change a bit our viewpoint on the definition of these random fields $H_N$ and $H_N^\mathrm{bip}$. Instead of writing them down explicitly as in \eqref{e.def.HN} and \eqref{e.def.HN.bip}, an equivalent way to define them is to specify that they are centered Gaussian fields, and to display their covariance. For the SK model, we have for every $\sigma, \tau \in \{\pm 1\}^N$ that
\begin{equation}  
\label{e.cov.HN}
\E \Ll[ H_N(\sigma) H_N(\tau) \Rr] = N \Ll( \frac{\sigma \cdot \tau}{N} \Rr) ^2,
\end{equation}
where the dot in $\sigma \cdot \tau$ denotes the usual scalar product in $\R^N$. So instead of writing the formula in \eqref{e.def.HN}, we could also have said ``let $(H_N(\sigma))_{\sigma \in \{\pm 1\}^N}$ be the centered Gaussian vector whose covariance is given by \eqref{e.cov.HN}''. The one advantage of the explicit formula, besides its possibly more intuitive appeal, is that it makes it transparent that such a random vector exists. More generally, one could consider centered Gaussian fields $(H_N(\sigma))_{\sigma \in \R^N}$ such that, for some smooth function $\xi : \R \to \R$, we have for every $\sigma, \tau \in \R^N$ that
\begin{equation}  
\label{e.def.cov}
\E \Ll[ H_N(\sigma) H_N(\tau) \Rr] = N \xi\Ll( \frac{\sigma \cdot \tau}{N} \Rr);
\end{equation}
the SK model corresponds to the case when $\xi(r) = r^2$. For the choice of $\xi(r) = r^3$, we can find a Gaussian field that satisfies \eqref{e.def.cov} by setting
\begin{equation*}  
H_N(\sigma) := \frac 1 N \sum_{i,j,k = 1}^N J_{i,j,k} \sigma_i \sigma_j \sigma_k,
\end{equation*}
where $(J_{i,j,k})$ are independent centered Gaussians with unit variance. One can realize $\xi(r) = r^p$ for any positive integer $p$ by proceeding similarly. Multiplying a random field by a factor of $\lambda$ transforms the function $\xi$ into $\lambda^2 \xi$. And by adding independent versions of fields with possibly different functions $\xi$, we can also create a new field whose covariance is the sum of those functions $\xi$. So we see that any function $\xi$ that can be written in the form 
\begin{equation}  
\label{e.valid.xi}
\xi(r) = \sum_{p =0}^{+\infty} a_p r^p,
\end{equation}
with $a_p \ge 0$ going to zero sufficiently rapidly as $p$ tends to infinity, will be a valid covariance function. One can show that these are all the admissible functions (see \cite[Proposition~6.6]{mourrat2023free} for a statement that also covers the case of models with multiple types of spins).

 For the bipartite model, we have instead that, for every $\sigma, \tau \in \{\pm 1\}^N \times \{\pm 1\}^N$,
\begin{equation}  
\label{e.def.cov.bip}
\E \Ll[ H_N^\bip(\sigma) H_N^\bip(\tau) \Rr] = N \Ll( \frac{\sigma_1 \cdot \tau_1}{N} \Rr) \Ll( \frac{\sigma_2 \cdot \tau_2}{N} \Rr).
\end{equation}
The key technical difference between the SK and the bipartite models is that here the relevant function that shows up on the right side of \eqref{e.def.cov.bip} is the mapping $(x,y) \mapsto xy$, which is \emph{not convex}. To be precise, for models with only one type of spins, i.e.\ of the form in \eqref{e.def.cov}, what is crucial is that the function $\xi$ is convex over $\R_+$; as one can see from~\eqref{e.valid.xi}, this is in fact always true! This convexity property can however break down as soon as we consider models with two or more types of spins. In general, we may consider models with a fixed number $D$ of types of spins, say $\sigma = (\sigma_1, \ldots, \sigma_D) \in (\R^N)^D$, with a covariance such that, for every $\sigma, \tau \in (\R^N)^D$,
\begin{equation*}  
\E \Ll[ H_N(\sigma) H_N(\tau) \Rr] = N \xi \Ll( \Ll( \frac{\sigma_d \cdot \sigma_{d'}}{N} \Rr)_{1 \le d, d' \le D}  \Rr) ,
\end{equation*}
where $\xi$ is some (admissible) function from $\R^{D\times D}$ to $\R$. Those models for which we can write down and rigorously prove a Parisi formula for the limit free energy are exactly those for which the function $\xi$ is convex over the space of symmetric positive semidefinite matrices (see \cite[Theorem~1]{chen2024free}, which crucially builds upon \cite{barcon, pan.aom, pan.multi, pan.potts, pan.vec}).

%
%
%
%
%
%
\section{A connection with partial differential equations}
\label{s.pdes}

So what can we do for the bipartite model? As discussed in the previous section, in this case it is not even clear what one is supposed to show, as there is no clear way to extend the Parisi formula into some reasonable candidate for the limit free energy. I will explain one possible route which my collaborators and I have been exploring. The idea is in two steps. First, we enrich the free energy, adding some (hopefully not too complicated) terms to the energy function, so that we end up with a free energy that depends on additional parameters besides the inverse temperature. Next, one hopes to find a partial differential equation that this free energy solves approximately, and to characterize the limit free energy as the unique solution to this equation. This approach works well for simpler models such as the Curie-Weiss model and its generalizations, or for some problems of statistical inference such as community detection, see \cite[Chapters 1 to 4]{HJbook} for a detailed presentation. Perhaps good pedagogic practice would require that I discuss these models first. Instead I will try to directly address the more complicated case of spin glasses, but some aspects of the discussion will then be rather sketchy. 

The fact that the limit free energy of a model of statistical mechanics may solve a partial differential equation is an old observation going back at least to \cite{bra83, new86} (see also \cite{bauerschmidt2023stochastic} for a recent survey on related topics). In the context of spin glasses, these connections were first explored in \cite{abarra, barra2,barra1, guerra2001sum} under simplifying assumptions. The fact that the Parisi formula can be recast as the value of the solution to some partial differential equation is from \cite{chen2022hamilton, mourrat2022parisi, mourrat2020extending} (see also \cite{HJbook}). 

In order to keep the notation simple, we will first discuss the approach in the case of the SK model. Let us start by defining a new free energy with an additional parameter, and then discuss motivations for this choice. For every $t \ge 0$ and $h \ge 0$, we set
\begin{equation}  
\label{e.def.newFN}
F_N(t,h) := - \frac 1 N \E \log \sum_{\sigma \in \{\pm 1\}^N} \exp(\sqrt{2t} H_N(\sigma) - Nt + \sqrt{2h} z \cdot \sigma - N h), 
\end{equation}
where $z = (z_1,\ldots, z_N)$ is a vector of independent centered Gaussians with unit variance, independent of $H_N$, and we recall that the function $H_N$ for the SK model is defined in~\eqref{e.def.HN}. 

Several comments need to be made to explain some of the choices that have been made in this new definition of $F_N$.  First, we replaced $\beta$ with $\sqrt{2t}$. The ultimate reason for this is that the subsequent formulas will look nicer in this way, but one way to sense that this may be so is to recall that $H_N(\sigma)$ is a Gaussian random variable, and so the scaling $\sqrt{2t} H_N(\sigma)$ is such that the variance of the Gaussian scales linearly, as with Brownian motion; it is as if we were continuously adding new independent copies of $H_N(\sigma)$ homogeneously in time. The same goes for $\sqrt{2h}$ in front of $z \cdot \sigma$. The factor $-Nt$ is also ultimately here just to make the subesquent formulas a bit nicer; it is half the variance of the Gaussian random variable $\sqrt{2t} H_N(\sigma)$. People who are familiar with stochastic calculus will recognize a stochastic exponential here: in particular, we have that $\E[\exp(\sqrt{2t} H_N(\sigma) - Nt)] = 1$. So if we were to exchange the expectation and the logarithm in \eqref{e.def.newFN}, we would end up with summands that are all equal to $1$. By Jensen's inequality, the quantity with expectation and logarithm interchanged is smaller than the original one, and we set things up so as to monitor the defect in this Jensen's inequality\footnote{In general, we also do not restrict ourselves to models defined on $\{\pm 1\}^N$. For a model whose covariance is given by \eqref{e.def.cov}, we thus write 
\begin{equation*}  
F_N(t,h) := - \frac 1 N \E \log \int \exp\Ll(\sqrt{2t} H_N(\sigma) - N t\xi(|\sigma|^2/N) + \sqrt{2h} z \cdot \sigma - h|\sigma|^2\Rr) \, \d P_N(\sigma),
\end{equation*}
where $P_N = P_1^{\otimes N}$ is the $N$-fold tensor product of a probability measure $P_1$ on $\R$ with compact support. We recover the SK model by choosing $\xi(r) = r^2$ and $P_1 = (\de_1 + \de_{-1})/2$, up to the addition of a trivial factor of $\log 2$. Notice that with this definition, Jensen's inequality tells us that $F_N \ge 0$, thanks to the extra minus sign that we added there. 
}.

From the formula \eqref{e.def.newFN} we have so far mostly discussed simple changes of variables, or tried to give some excuse for the presence of the silly-looking terms $Nt$ or $Nh$, but the important part is to get a feeling as to why we are adding this term $\sqrt{2h} z \cdot \sigma$. That we are looking for a partial differential equation for $F_N$ means that we hope to be able to compensate small variations of $t$ with small variations in $h$. So we aim to find some term that looks like $H_N(\sigma)$ in some sense, while being also simpler to analyse. Perhaps one way to think that the term $\sqrt{2h} z\cdot \sigma$ is not an unreasonable choice is to write 
\begin{equation*}  
H_N(\sigma) = \frac 1 {\sqrt N} \sum_{i = 1}^N \Ll( \sum_{j = 1}^N J_{ij} \sigma_j\Rr) \sigma_i ,
\end{equation*}
and to venture the guess that maybe the random variables $(\sum_{j = 1}^N J_{ij} \sigma_j)$ could be substituted with equivalent independent Gaussians, because, well, at least for each fixed $\sigma$ they are Gaussian after all. A more detailed explanation is beyond the scope of this note; for me the best heuristic is that discussed in \cite[Exercise~6.5 and solution]{HJbook}. 

Let us be pragmatic here and just calculate the derivatives of $F_N$ to see if something interesting happens. In order to express these derivatives nicely, we introduce some notation for the Gibbs measure. For any function $f$, we write
\begin{equation}  
\label{e.def.Gibbs}
\langle f(\sigma) \rangle := \frac{\sum_{\sigma \in \{\pm 1\}^N} f(\sigma) \exp(H_N(t,h,\sigma))}{\sum_{\sigma \in \{\pm 1\}^N} \exp(H_N(t,h,\sigma))},
\end{equation}
where we have set $H_N(t,h,\sigma) := \sqrt{2t} H_N(\sigma) - Nt + \sqrt{2h} z \cdot \sigma - N h$. In the notation on the left side of \eqref{e.def.Gibbs}, the bracket $\langle \cdot \rangle$ thus stands for the expectation with respect to the Gibbs measure, and we think of $\sigma$ as a random variable that is sampled accordingly. We also write $\sigma'$ to denote an independent copy of $\sigma$ under the Gibbs measure, that is, we write
\begin{equation*}  
\langle f(\sigma, \sigma') \rangle := \frac{\sum_{\sigma,\sigma' \in \{\pm 1\}^N} f(\sigma, \sigma') \exp(H_N(t,h,\sigma)+H_N(t,h,\sigma'))}{\sum_{\sigma,\sigma' \in \{\pm 1\}^N} \exp(H_N(t,h,\sigma)+H_N(t,h,\sigma'))}.
\end{equation*}
A calculation\footnote{You can try! The only thing to have in mind is that if $G$ is a centered Gaussian of unit variance and $f$ is a nice enough function (say $C^1$ with reasonable growth of $f$ and $f'$ at infinity), then 
\begin{equation*}  
\E[G f(G)] = \E[f'(G)].
\end{equation*}
To see this, you just need to write the expectation explicitly as an integral against the Gaussian density, and integrate by parts.} gives us that
\begin{equation}  
\label{e.FN.derivatives}
\dr_t F_N(t,h) = \E \la \Ll( \frac{\sigma \cdot \sigma'}{N} \Rr) ^2 \ra \quad \text{ and } \quad \dr_h F_N(t,h) = \E \la  \frac{\sigma \cdot \sigma'}{N}  \ra  . 
\end{equation}
The dependence on $t$ and $h$ in the right-hand sides of the identities above is hidden in the definition of the Gibbs average $\la \cdot \ra$.
For more general models as in \eqref{e.def.cov}, we would find the exact same expression for $\dr_h F_N$ as in \eqref{e.FN.derivatives} (for the corresponding definition of the Gibbs measure), while for the derivative in $t$, we would find that
\begin{equation}  
\label{e.drt.xi}
\dr_t F_N(t,h) = \E \la \xi \Ll( \frac{\sigma \cdot \sigma'}{N} \Rr)  \ra.
\end{equation}
Continuing with the SK model for now, we thus obtain that 
\begin{equation}
\label{e.pde.FN}
\dr_t F_N - (\dr_h F_N)^2 = \E \la \Ll( \frac{\sigma \cdot \sigma'}{N} \Rr) ^2 \ra  - \Ll( \E \la  \frac{\sigma \cdot \sigma'}{N}  \ra   \Rr) ^2.
\end{equation}
The right-hand side of this identity is the variance of the random variable $\sigma \cdot \sigma'/{N}$ under $\E \la \cdot \ra$. Since $\sigma \cdot \sigma'/N$ is a sum of a large number of terms, we could at first anticipate that it will have small fluctuations --- and in fact, this intuition turns out to be valid as long as $t$ is small. Just to see what happens (or by deciding that we restrict to small $t$), let us get along a little with this hypothesis that the variance of $\sigma \cdot \sigma'/{N}$ tends to zero as $N$ tends to infinity. We are then led to the expectation that $F_N$ may converge to a limit function~$f$ that solves the equation
\begin{equation}  
\label{e.hj.SK.simple}
\dr_t f - (\dr_h f)^2 = 0.
\end{equation}
Moreover, as we set $t = 0$ in the definition of $F_N$ in \eqref{e.def.newFN}, we find that we can easily calculate the result by writing $z \cdot \sigma = \sum_{i = 1}^N z_i \sigma_i$, writing the exponential of the sum as a product of exponentials, and realizing that we can then factorize the summation into a product of $N$ simple sums over $\{\pm 1\}$. In short, we see that 
\begin{equation}
\label{e.simple.init}
F_N(0,h) = F_1(0,h). 
\end{equation}
So if we believe that the random variable $\sigma \cdot \sigma'/N$ has vanishingly small fluctuations in the limit of large $N$, then we are tempted towards the conjecture that $F_N$ converges to the function $f$ that solves \eqref{e.hj.SK.simple} with initial condition $f(0,\cdot) = F_1(0,h)$. In fact, this guess for the limit free energy that we thus obtain is exactly the one that was proposed (using other arguments) in the paper that introduced the SK model \cite{sherrington1975solvable}. For the model with the covariance as in~\eqref{e.def.cov}, we would instead have guessed the partial differential equation
\begin{equation}
\label{e.hj.xi.simple}
\dr_t f  - \xi(\dr_h f) = 0.
\end{equation}
Partial differential equations of this sort are called Hamilton-Jacobi equations.

The hypothesis that the random variable $\sigma \cdot \sigma'/N$ has vanishingly small fluctuations in the limit of large $N$ is just of the sort that works very well for simpler models, but here it is only valid at high temperature, that is for small values of $t$. To see why, notice first that since $|\sigma| = |\sigma'| = \sqrt{N}$, the fluctuations of $\sigma \cdot \sigma'/N$ are the same as those of $|\sigma-\sigma'|^2/(2N)$. The point then is that for sufficiently large values of $t$, the Gibbs measure starts to look like the ultrametric tree appearing on Figure~\ref{f.ultrametric}, with each leaf vertex endowed with a fixed probability weight. As long as this tree is not reduced to a single vertex (which does happen for large $t$), we see that the distance between two points taken at random on the leaves of the tree will fluctuate. 

Even though our first attempt at analyzing the model did not really work out, the idea we have explored here seems promising, so we are going to build upon it. The point is that we have not fully succeeded in ``closing the equation'' with the introduction of this parameter $h$; but maybe if we were introducing more parameters, then we would be able to do better. 

The full construction of the free energy that incorporates these additional parameters would take too much space to explain here (one can consult \cite[Section~6.4]{HJbook}), but the basic idea is that, since we anticipate that the true system may have a complicated ultrametric structure as discussed in the previous section, we want to replace our naive ``external field'' $\sqrt{2h} z$, which we were ``dotting'' against $\sigma$, with a more refined object that has such an ultrametric structure embedded into it. This ultrametric structure is encoded using a non-decreasing function $q : [0,1] \to [0,\infty)$.

So this construction leads us to a new definition of the free energy $F_N$, which we now think of as a function of $t$ and $q$ in place of $t$ and $h$. We may call this extended function $F_N(t,q)$ the ``enriched'' free energy. The version in \eqref{e.def.newFN} corresponds to the choice of the constant function $q = h$ in this new definition; in particular, we are computing the free energy of the ``vanilla'' model when we compute $F_N(\beta^2/2,0)$, up to remembering the extra minus sign and the term $Nt$ in \eqref{e.def.newFN}. 

Let us denote by $\mcl Q$ the space of non-decreasing functions $q : [0,1] \to [0,\infty)$. The notion of differentiation for functions $g : \mcl Q \to \R$ that we use is as follows. We say that $g$ is differentiable at $q \in \mcl Q$, and in this case denote by $\dr_q g(q,\cdot) \in L^2([0,1];\R)$ the derivative at $q$, if for every $q' \in \mcl Q$, we have 
\begin{equation*}  
g((1-\ep) q + \ep q') = g(q) + \ep \int_0^1 (q'-q)(u) \dr_q g(q,u) \, \d u + o(\eps) \qquad (\ep \to 0). 
\end{equation*}
Equipped with this, one can find a nice expression for the derivative of $F_N$ with respect to $q$, which resembles what we have found in the second part of \eqref{e.FN.derivatives} but also involves some other variables that play a role in the ultrametric structure of the new external field. One is then led to form the quantity
\begin{equation*}  
\dr_t F_N(t,q) - \int_0^1 (\dr_q F_N(t,q,u))^2 \, \d u = \text{small}?
\end{equation*}
What replaces the ``small'' right-hand side above is a conditional variance of $\sigma \cdot \sigma'/N$, in particular it is indeed smaller than the variance of $\sigma \cdot \sigma'/N$, so we are indeed making progress. For the model with the covariance as in \eqref{e.def.cov}, we are led to hope that 
\begin{equation*}  
\dr_t F_N(t,q) - \int_0^1 \xi(\dr_q F_N(t,q,u)) \, \d u = \text{small}?
\end{equation*}
While we did not explain the construction of this new external field with an ultrametric structure, it turns out that it is still ``factorizable'' in the sense that the nice property that we had in \eqref{e.simple.init} is still valid in this more general setting, that is, we have for every $q \in \mcl Q$ that
\begin{equation}
\label{e.init}
F_N(0,q) = F_1(0,q).
\end{equation}
It turns out that our more sophisticated hopes are now correct. For notational convenience, we define $\psi(q) := F_1(0,q)$. 
\begin{theorem}[The Parisi formula as a Hamilton-Jacobi equation \cite{chen2022hamilton, mourrat2022parisi, mourrat2020extending}]
\label{t.hj}
The enriched free energy $F_N : \R_+ \times \mcl Q \to \R$ for the SK model converges pointwise to the function $f : \R_+ \times \mcl Q \to \R$ that solves 
\begin{equation}
\label{e.hj.SK}
\begin{cases}
\dr_t f - \int_0^1 (\dr_q f)^2  = 0   & \text{ on } \R_+ \times \mcl Q, \\
f(0,\cdot) = \psi  & \text{ on } \mcl Q.
\end{cases}
\end{equation}
More generally, if $F_N : \R_+ \times \mcl Q \to \R$ now stands for the enriched free energy associated with a model whose covariance is given by \eqref{e.def.cov}, then $F_N$ converges to the function $f$ that solves
\begin{equation}
\label{e.hj.xi}
\begin{cases}
\dr_t f - \int_0^1 \xi(\dr_q f)  = 0   & \text{ on } \R_+ \times \mcl Q, \\
f(0,\cdot) = \psi  & \text{ on } \mcl Q.
\end{cases}
\end{equation}
\end{theorem}
Part of the task of making sense of this theorem is that one needs to find a good notion of solution for the Hamilton-Jacobi equations in \eqref{e.hj.SK} and \eqref{e.hj.xi}. This is based on the notion of viscosity solutions (see \cite[Chapter 3]{HJbook} for an introduction that is tailored to our context). 

How does this formulation relate to the Parisi formula? The key point is that when the non-linearity (the square function in \eqref{e.hj.SK}) is convex, the solution admits a variational representation known as the Hopf-Lax formula. The solution to \eqref{e.hj.SK} can be written as
\begin{equation}  
\label{e.hopf-lax.SK}
f(t,q) = \sup_{q' \in \mcl Q} \Ll(\psi(q+q') - \frac 1 {4t} \int_0^1 (q')^2 \Rr).
\end{equation}
For the solution to the more general equation in \eqref{e.hj.xi}, we have
\begin{equation}  
\label{e.hopf-lax.xi}
f(t,q) = \sup_{q' \in \mcl Q} \Ll(\psi(q+q') - t \int_0^1 \xi^* \Ll( \frac{q'}{t} \Rr)  \Rr),
\end{equation}
where $\xi^*$ is the convex dual of $\xi$ defined, for every $s \in \R$, by
\begin{equation*}  
\xi^*(s) := \sup_{r \ge 0} \Ll( rs - \xi(r) \Rr) .
\end{equation*}
One can then bridge this representation back to the Parisi formula by setting $q = 0$ in~\eqref{e.hopf-lax.SK}, making a change of variables, and doing some explicit calculations involving the function~$\psi$.\footnote{Our formulas here are with a supremum, unlike in \eqref{e.parisi.simple}, but this is only for the trivial reason that we have added a minus sign in the definition in \eqref{e.def.newFN}.}

All this is nice, but we already knew the answer in these cases. The really interesting things come about when we consider models with different types of spins such as the bipartite model. Recall that in this case the energy function $H_N^\bip$ is defined in~\eqref{e.def.HN.bip}. Let us first see how a simpler guess analogous to \eqref{e.hj.SK.simple} or \eqref{e.hj.xi.simple} would look like then. The important point is that since there are now two types of spins, we would like to make sure that we have one additional variable (one additional ``external field'') to act upon each of these two types. So we set, for every $t \ge 0$, $h = (h_1,h_2) \in \R_+^2$, and $\sigma = (\sigma_1, \sigma_2) \in \{\pm 1\}^N \times \{\pm 1\}^N$,
\begin{equation*}  
H_N(t,h,\sigma) := \sqrt{2t} H_N^\bip(\sigma) - Nt + \sqrt{2h_1} z_1 \cdot \sigma_1 - N h_1 + \sqrt{2 h_2} z_2 \cdot \sigma_2 - N h_2,
\end{equation*}
(we drop the superscript $^\bip$ from now on for ease of notation), and then for every $t \ge 0$ and $h = (h_1,h_2) \in \R_+^2$, we set the new free energy to be
\begin{equation*}  
F_N(t,h) := - \frac 1 N \E \log \sum_{\sigma \in \{\pm 1\}^N \times \{\pm 1\}^N} \exp(H_N(t,h,\sigma)).
\end{equation*}
The definition of the Gibbs average $\la \cdot \ra$ is hopefully clear, it is essentially as in the formula in \eqref{e.def.Gibbs}, except that now the variable $h$ is in $\R_+^2$ and the summation variable $\sigma$ ranges in $\{\pm 1\}^N \times \{\pm 1\}^N$. In place of \eqref{e.FN.derivatives} or \eqref{e.drt.xi}, we obtain that 
\begin{equation*}  
\dr_t F_N(t,h) = \E \la \Ll( \frac{\sigma_1 \cdot \sigma_1'}{N} \Rr) \Ll( \frac{\sigma_2 \cdot \sigma_2'}{N} \Rr) \ra,
\end{equation*}
with still
\begin{equation*}  
\dr_{h_1} F_N(t,h) = \E \la  \frac{\sigma_1 \cdot \sigma_1'}{N}  \ra  \quad \text{ and } \quad \dr_{h_2} F_N(t,h) = \E \la  \frac{\sigma_2 \cdot \sigma_2'}{N}  \ra  .
\end{equation*}
So now if we expect that the random variables $\sigma_1 \cdot \sigma_1'/N$ and $\sigma_2 \cdot \sigma_2'/N$ do not fluctuate much in the limit of large $N$, we are led to the belief that $F_N$ converges to $f$ solving 
\begin{equation*}  
\dr_t f - \dr_{h_1} f \,  \dr_{h_2} f = 0. 
\end{equation*}
The initial condition $f(0,\cdot)$ is easy to compute as it satisfies \eqref{e.simple.init} again. 

The problem is that this guess is too naive again of course, so instead we need to replace our pair of simple variables $h = (h_1, h_2) \in \R_+^2$ by a pair of functions $q = (q_1, q_2) \in \mcl Q^2$ that encodes the ultrametric structure of the external fields that act upon $\sigma_1$ and $\sigma_2$ respectively. This defines for us a more sophisticated free energy $F_N$, now defined on $\R_+ \times \mcl Q^2$ in place of $\R_+ \times \R_+^2$. Again we have that \eqref{e.init} holds (now with $q$ ranging in $\mcl Q^2$); let us set $\psi := F_1(0,\cdot)$ for convenience. By going through similar calculations as for the SK model, one thus ends up expecting that:
\begin{conjecture}  
\label{conj}
Maybe the enriched free energy $F_N : \R_+ \times \mcl Q^2 \to \R$ for the bipartite model converges to the function $f : \R_+ \times \mcl Q^2 \to \R$ that solves
\begin{equation}
\label{e.hj.bip}
\begin{cases}
\dr_t f - \int_0^1 \dr_{q_1} f \, \dr_{q_2} f  = 0   & \text{ on } \R_+ \times \mcl Q^2, \\
f(0,\cdot) = \psi  & \text{ on } \mcl Q^2 \qquad ?
\end{cases}
\end{equation}
\end{conjecture}
In this case, since the non-linearity in the equation \eqref{e.hj.bip} (i.e. the mapping $(x,y) \mapsto xy$) is neither convex nor concave, there is no standard way to represent the solution to \eqref{e.hj.bip} variationally. So there is no clear way to rewrite this candidate for the limit free energy into a form that would resemble the Parisi formula for the SK model. 

\begin{figure}
\centering
\includegraphics[width=0.36\linewidth]{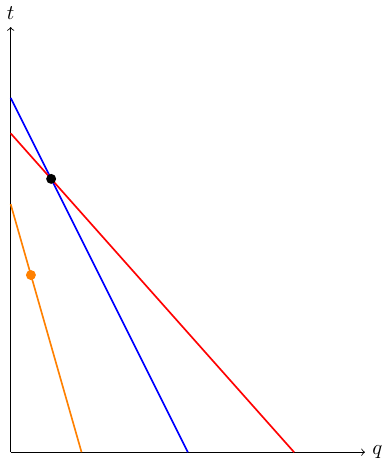}
\captionsetup{width=0.95\linewidth}
\caption{{\small Each characteristic line offers us a ``prediction'' for what the limit free energy is. For small $t$, each point $(t,q)$  (e.g. the orange point) is reached by exactly one characteristic line, so we know the value of the free energy then. For large $t$, it could happen that multiple characteristics reach a point, as happens for the black point. We know that the limit free energy must be as prescribed by one of those characteristics, but we do not know which one.}}
    \label{f.charac}
\end{figure}

As was said already, we do not know if Conjecture~\ref{conj} is valid or not. One thing we know is that any subsequential limit of $F_N$ must satisfy the equation in~\eqref{e.hj.bip} ``almost everywhere'' in $\R_+ \times \mcl Q^2$, for a suitable interpretation of ``almost everywhere'' \cite{chen2024free}. This is not sufficient to fully characterize the solution to \eqref{e.hj.bip} though; even in finite dimensions, there may be many functions that solve a Hamilton-Jacobi equation almost everywhere. Another thing that we know is that the solution to \eqref{e.hj.bip} is a lower bound \cite{chen2022hamilton, mourrat2020nonconvex}, that is, for every $t \ge 0$ and $q \in \mcl Q^2$, we have
\begin{equation*}  
\liminf_{N \to \infty} F_N(t,q) \ge f(t,q),
\end{equation*}
where $f$ denotes the unique viscosity solution to \eqref{e.hj.bip}. The result that I find most interesting, and which essentially matches what the physicists say using other language, is one that I will only describe informally. In order to do so, let us first discuss a classical method to solve \eqref{e.hj.bip} for short times called the method of characteristics (see \cite[Section~3.5]{HJbook} for more precision). It turns out that if one does some formal calculations assuming that the solution to \eqref{e.hj.bip} is twice differentiable, one can guess a simple formula for what the value of the solution ought to be along each of a spanning family of lines called the characteristics. In our case the characteristic that starts at $q$ is the straight line
\begin{equation*}  
t \mapsto q - t\nabla \xi(\dr_q \psi(q)),
\end{equation*}
where we denoted by $\nabla \xi : \R_+^2 \to \R_+^2$ the function mapping $(x,y)$ to $(y,x)$ (in line with the idea that the relevant ``function $\xi$'' for the bipartite model is $(x,y) \mapsto xy$), and $\dr_q \psi = (\dr_{q_1} \psi, \dr_{q_2} \psi)$. 
This guess provided by characteristics is actually really valid for short times, in the sense that it really gives us the value of the viscosity solution there. The problem is that for large times, these lines may start to intersect each other, and they do not agree on what they predict (and the viscosity solution may then take a value that is not even among those multiple options). One of the main results of \cite{chen2024free} is that if one assumes that the enriched free energy $F_N : \R_+ \times \mcl Q^2 \to \R$ of the bipartite model converges to some limit function $f$, then there is always at least one characteristic line that prescribes the correct value for $f$. The only remaining source of ambiguity is that we are not able to say which of the characteristic lines is the correct one, in case multiple lines arrive at a given point. See also Figure~\ref{f.charac} for an illustration.

%
%
%
%
%
%
\section{Closing thoughts}

Spin glasses are models of statistical mechanics that are relatively simple to define. Yet, one aspect I find striking about them is that they display a very rich and interesting mathematical structure. Despite their apparent simplicity, they may give us some insight into the behavior of various ``complex'' systems across disciplines. 

I hope that this informal overview of the Parisi formula and of some related open questions has piqued your interest. I want to stress though that this is only one facet of the study of spin glasses, and there surely is something for everyone there! Books on spin glasses include \cite{bolthausen2007spin, bovier2006statistical, bovier1998mathematical, charbonneau2023spin, contucci2013perspectives, dominicis2006random, HJbook, mpvbook, nishimori2001statistical, opper2001advanced, pan, stein2013spin, Tbook1, Tbook2}. The book~\cite{HJbook} puts particular emphasis on the PDE point of view discussed here in Section~\ref{s.pdes} (and contains an overview of the research literature). For a multi-disciplinary perspective and applications of spin glasses to other areas, one can consult the recent book \cite{charbonneau2023spin} on ``spin glass theory and far beyond''.

\small
\bibliographystyle{plain}
\bibliography{gazette}

\end{document}